\documentclass[a4paper,leqno,11pt]{amsart}
\usepackage{amsfonts,amssymb,verbatim,amsmath,amsthm,latexsym,textcomp,amscd}
\usepackage{latexsym,amsfonts,amssymb,epsfig,verbatim}
\usepackage{amsmath,amsthm,amssymb,latexsym,graphics,textcomp}
\usepackage{graphicx}
\usepackage{color}
\usepackage{url}
\usepackage{enumerate}
\usepackage[mathscr]{euscript}
\usepackage{tikz}
\usepackage{dsfont}
\usetikzlibrary{matrix}
\usepackage{sseq}
\usepackage{hyperref}

\input xy
\xyoption{all}

\theoremstyle{plain}
\newtheorem*{thma}{Theorem A}
\newtheorem*{thmb}{Theorem B}

\theoremstyle{definition}
\newtheorem{theorem}{Theorem}[section]
\newtheorem{prop}[theorem]{Proposition}
\newtheorem{lemma}[theorem]{Lemma}
\newtheorem{cor}[theorem]{Corollary}

\newtheorem{subsec}[theorem]{}

\newtheorem{exam}[theorem]{Example}

\newtheorem{thm}[theorem]{Theorem}

\theoremstyle{remark}
\swapnumbers
   \newtheorem{ack}[theorem]{Acknowledgements}

\newcommand{\C}{{\mathbb C}}

\newcommand{\Z}{{\mathbb{Z}}}

\newcommand{\Q}{{\mathbb Q}}

\newcommand\FF{{\mathcal F}}

\newcommand\LL{{\mathcal L}}
\newcommand\MM{{\mathcal M}}

\newcommand\PP{{\mathcal P}}

\newcommand\PMF{{\PP\kern-2pt\MM\FF}}

\newcommand\PML{{\PP\kern-2pt\MM\LL}}

\newcommand{\fsubd}{\mathrel{{\scriptstyle\searrow}\kern-1ex^d\kern0.5ex}}
\newcommand{\bsubd}{\mathrel{{\scriptstyle\swarrow}\kern-1.6ex^d\kern0.8ex}}
\newcommand{\fsubeq}{\mathrel{\raise-.7ex\hbox{$\overset{\searrow}{=}$}}}
\newcommand{\bsubeq}{\mathrel{\raise-.7ex\hbox{$\overset{\swarrow}{=}$}}}

\newcommand{\tsh}[1]{\left\{\kern-.9ex\left\{#1\right\}\kern-.9ex\right\}}

\newenvironment{myeq}[1][]
{\stepcounter{theorem}\begin{equation}\tag{\thetheorem}{#1}}
{\end{equation}}

\newenvironment{mysubsection}[2][]
{\begin{subsec}\begin{upshape}\begin{bfseries}{#2.}
			\end{bfseries}{#1}}
		{\end{upshape}\end{subsec}}

\makeatletter
\@namedef{subjclassname@2020}{%
  \textup{2020} Mathematics Subject Classification}
\makeatother

\title{$BP$-Cohomology of Projective Stiefel Manifolds}
\author{Samik Basu, Debanil Dasgupta}

\email{samik.basu2@gmail.com; samikbasu@isical.ac.in}
\address{Stat-Math Unit,
Indian Statistical Institute,
B. T. Road, Kolkata-700108, India.}

\email{debanil12@gmail.com;}
\address{Stat-Math Unit,
Indian Statistical Institute,
B. T. Road, Kolkata-700108, India.}


\begin{document}

\maketitle

\begin{abstract}
In this paper, we compute the $BP$-cohomology of complex projective Stiefel manifolds. The method involves the homotopy fixed point spectral sequence, and works for complex oriented cohomology theories. We also use these calculations and $BP$-operations to prove new results about  equivariant maps between Stiefel manifolds.
\end{abstract}

\section{Introduction} 
The projective Stiefel manifolds have been of interest in connection with a varied spectrum of topological questions. On one hand, they are useful in studying equivariant maps between the Stiefel manifolds, and on the other, they form a part of an obstruction theory for constructing sections of multiples of a given line bundle. In the real case, they play an important role in the immersion problem for real projective spaces. 

In this paper, we consider the complex projective Stiefel manifolds $PW_{n,k}$, defined as the quotient of the complex Stiefel manifold $W_{n,k}$ by the $S^1$-action. The cohomology of $PW_{n,k}$ with $\Z/p$-coefficients was computed in \cite{AGMP}, which is analogous to the $\Z/2$-computation for real projective Stiefel manifolds in \cite{GH68}. Among other applications, this has been used to prove the non-existence of $S^1$-equivariant maps between Stiefel manifolds \cite{PP13}. 

A natural idea here is that extending the computations to generalized cohomology theories would yield further results about equivariant maps. We follow through along these lines and compute the $BP$-cohomology as (Theorem \ref{projstbp})
\begin{thma}
The $BP$-cohomology of $PW_{n,k}$ is described as
 $$ BP^*(PW_{n,k}) \cong \Lambda_{BP^*(pt)}(\gamma_{n-k+2},\cdots,\gamma_n)\otimes_{BP^*(pt)} BP^*(pt)[[x]]/I $$
 where $\gamma_j$'s are of degree $2j-1$, $x$ is of degree $2$, and $I$ is the ideal generated by  $\{\binom{n}{j}x^j\vert n-k<j\leq n \}$. 
\end{thma}
The method used to compute the $BP$-cohomology is the homotopy fixed point spectral sequence. This works for any complex oriented cohomology theory, where the class $x$ comes from the choice of complex orientation. Consequently, the $K$-theory of the complex projective Stiefel manifold has an analogous formula, which was computed in \cite{Gon20} using the Hodgkin spectral sequence for the cohomology of homogeneous spaces. The same method is also likely to work for $P_\ell W_{n,k}$, the quotient by a variant of the $S^1$-action, whose cohomology was computed in \cite{BS17}.   

We observe that the $BP$-cohomology ring of $PW_{n,k}$ is just the extension of coefficients from $\Z_{(p)}$ in ordinary cohomology to $\Z_{(p)}[v_1,v_2,\cdots]$ in $BP$-cohomology. Therefore, the primary multiplicative structure does not yield new results for equivariant maps between Stiefel manifolds. However, $BP$ has the action of Adams operations \cite{Ara75}, which yield the following new result on equivariant maps. (see Theorem \ref{eqstbp})
\begin{thmb}
Suppose that $m,n,l, k$ are positive integers satisfying\\ 
1) $n-k<m-l$ and there is an $s$ such that $m<2^s+m-l\leq n$. \\ 
2) $2$ divides all the binomial coefficients $\binom{n}{n-k+1},\cdots, \binom{n}{m-l}$. \\
3) $2$ does not divide $\binom{m}{m-l+1}$ and $2\nmid m-l$. \\
Then, there is no $S^1$-equivariant map from $W_{n,k}$ to $W_{m,l}$.
\end{thmb}
We also obtain some new results using the action of Steenrod operations on $H^\ast PW_{n,k}$. We point out that the analysis of equivariant maps on Stiefel manifolds also leads to results in topological combinatorics \cite{BK21}. 

\begin{mysubsection}{Organization}
In section \ref{hfp}, we discuss the construction of the homotopy fixed point spectral sequence, proving results about the convergence and the differentials in the case of projective Stiefel manifolds. In section \ref{cohwnk}, we describe the cohomology of the Stiefel manifold over generalized cohomology theories, specializing to $BP$. In Section \ref{cohpwnk}, we complete the calculation of the spectral sequence and describe the $BP$-cohomology of $PW_{n,k}$. In section \ref{eqstmap}, we discuss the applications to equivariant maps between Stiefel manifolds. 
\end{mysubsection}

\begin{ack}
The first author would like to thank Shilpa Gondhali for some helpful conversations during the initial part of this work. The second author would like to thank Aloke Kumar Ghosh. The research of the first author was supported by the SERB MATRICS grant 2018/000845.
\end{ack}

\section{Homotopy fixed point spectral sequence}\label{hfp}

The purpose of this section is to set up the computational tools for the following sections. The main idea here is the homotopy fixed point spectral sequence for (naive) $G$-equivariant spectra: for a spectrum $Z$ with a $G$-action, there is a spectral sequence with $E_2$-page $H^s(G;\pi_{t} Z) $ which converges to $\pi_{t-s} Z^{hG}$ \cite{Dug03}. 

The principal example for our paper is when $G$ acts on a function spectrum $F(X,E)$ for a spectrum $E$\footnote{A spectrum here refers to one with trivial $G$-action.} and a based $G$-space $X$\footnote{This means that the base-point is fixed by $G$.}. Let us make this more precise. Let $E$ be a spectrum so that the reduced $E$-cohomology of based spaces are computed as 
\[
\tilde{E}^n(X) \cong [X, \Sigma^n E] \cong \pi_{-n}F(X,E).
\]
Here we use the notation $[-,-]$ for the homotopy classes of maps between spectra. We follow the construction of function spectra in \cite{May96}. If $X$ has a $G$-action, the function spectrum $F(X,E)$ is a spectrum with $G$-action (that is, a $G$-spectrum indexed over a trivial $G$-universe). We write $[-,-]^G_{tr}$ for the equivariant homotopy classes in the category of spectra with $G$-action, and $F^G_{tr}(-,-)$ for the equivariant function spectrum with $G$-action. We have the following result regarding this construction. 
\begin{prop}\label{proj-Stief-red} \cite[Ch. XVI, \S 1, (1.9)]{May96}
Let $X$ be a based $G$-space, and $E$ a spectrum. Then, 
\[
\pi_{-n}^G F^G_{tr}(X, E) \cong [X,\Sigma^n E]^G_{tr} \cong [X/G, \Sigma^n E] \cong \tilde{E}^n(X/G).
\]
\end{prop} 

For a free $G$-space $X$ we may apply Proposition \ref{proj-Stief-red} by adding a disjoint base-point. The homotopy fixed points of a spectrum $Z$ with $G$-action are $Z^{hG} = F^G_{tr}(EG_+,Z)^G$. We know that for a free $G$-space $X$, the projection $X\times EG\to X$ is a $G$-equivalence. Therefore, we have the following equivalence of spectra. 
\begin{cor} \label{free-sp-red} 
Let $X$ be a free $G$-space, and $E$ a spectrum. Then 
\[ 
F^G_{tr}(X_+,E)^{hG} \simeq F^G_{tr}(X_+,E)^G \simeq F(X/G_+, E).  
\]
\end{cor}

In this paper, we apply Corollary \ref{free-sp-red} to the case $X=W_{n,k}$, the Stiefel manifold of $k$-orthogonal vectors in $\C^n$. The group $G=S^1$ which acts on $W_{n,k}$ by vector-wise multiplication using $S^1 \subset \C$. This action is free and the quotient space is the projective Stiefel manifold $PW_{n,k}$. 
\begin{cor} \label{proj-stief-htpyfix} 
Let $E$ be a spectrum. There is an equivalence of spectra 
\[
F({PW_{n,k}}_+,E) \simeq F^{S^1}_{tr}({W_{n,k}}_+, E)^{hS^1}. 
\]
\end{cor}

We attempt to understand $F^{S^1}_{tr}({W_{n,k}}_+, E)$ via the homotopy fixed point spectral sequence. For a spectrum $Z$ with $S^1$-action, we follow the exposition in \cite{BR05} replacing homology with homotopy groups. We have a $S^1$-equivariant filtration of $ES^1$ given by 
\[
 \varnothing \subset S(\C) \subset S(\C^2) \subset \cdots S(\C^r) \subset S(\C^{r+1}) \subset \cdots, 
\] 
so that 
\[
Z^{hS^1} \simeq \varprojlim_r F^{S^1}_{tr}(S(\C^r)_+, Z)^{S^1}.
\]
We index the filtration of $ES^1$ as 
\[ 
E^{(r)}S^1 = \begin{cases}  
 S(\C^{\frac{r}{2}+1})  & \mbox{if } r \mbox{ is even} \\ 
 E^{(r-1)}S^1 & \mbox{if } r \mbox{ is odd}, 
 \end{cases} 
 \]
 so that 
 \[
 E^{(2r)}S^1/E^{(2r-1)} S^1 \simeq S^1_+ \wedge S^{2r}, ~~ E^{(2r+1)}S^1/E^{(2r)}S^1 \simeq \ast,
 \]
where the action of $S^1$ on $S^{2r}$ is the trivial action. The filtration on the induced tower of fibrations is written as 
 \[
 Z^{hS^1}_{(r)} =    F^{S^1}_{tr}(E^{(r)}S^1_+, Z)^{S^1},
 \]
 so that 
 \[                                      
 Z^{hS^1}_{(r)}/ Z^{hS^1}_{(r-1)} \simeq \begin{cases}  
 F^{S^1}_{tr}(S^1_+\wedge S^r, Z)^{S^1} \simeq \Sigma^{-r} Z   & \mbox{if } r \mbox{ is even} \\ 
 \ast  & \mbox{if } r \mbox{ is odd}. 
 \end{cases}
 \]
We may now follow \cite{BR05} to obtain a conditionally convergent spectral sequence \cite{Bo99}. 
\begin{prop}\label{htpyfpsseq}
Let $Z$ be a homotopy commutative ring spectrum with $S^1$-action. There is a conditionally convergent multiplicative spectral sequence 
\[
E_2^{s,t} = H^s( S^1; \pi_t(Z)) \implies \pi_{t-s}(Z^{hS^1}). 
\]
In this expression, the group cohomology $H^\ast(S^1;\pi_t Z)$ of $S^1$ in the discrete group $\pi_t Z$ equals $\Z[y] \otimes \pi_t Z$ with $|y|=(2,0)$. 
\end{prop}

\begin{exam} \label{z=e} 
If $Z=E$ with trivial $S^1$-action,  the homotopy fixed point spectrum $Z^{hS^1} \simeq F(BS^1_+, E)$. In this case, the homotopy fixed point spectral sequence becomes 
\[
E_2^{s,t} = H^s(\C P^\infty ) \otimes \pi_t E \implies \pi_{t-s} F(\C P^\infty_+, E).  
\]
Making identifications $E^n(\C P^\infty) \cong \pi_{-n}F(\C P^\infty_+,E)$, we observe that this reduces to the Atiyah-Hirzebruch spectral sequence for $\C P^\infty$. If $E$ is complex orientable, the class $y$ becomes a permanent cycle. 
\end{exam}

Next we specialize to the case $Z=F^{S^1}_{tr}(X_+, E)$ where $X$ is a free $S^1$-space, and $E$ is a spectrum. The homotopy groups of $F^{S^1}_{tr}(X_+,E)$ in Proposition \ref{htpyfpsseq} is computed by forgetting the $S^1$ action, and thus we have, 
\[
\pi_t F^{S^1}_{tr}(X_+, E) \cong \pi_t F(X_+,E) \cong E^{-t}(X).  
\]
On the other hand, we apply Corollary \ref{free-sp-red} to deduce
\[
\pi_t F^{S^1}_{tr}(X_+, E)^{hS^1} \cong \pi_t F^{S^1}_{tr}(X_+, E)^{S^1} \cong \pi_t F(X/S^1_+,E) \cong E^{-t}(X/S^1).  
\]
We now switch the sign of the $t$-grading in the spectral sequence of Proposition \ref{htpyfpsseq} to obtain a conditionally convergent multiplicative spectral sequence 
\[
E_2^{s,t} = H^s( S^1; E^t(X)) \cong \Z[y] \otimes E^t(X) \implies E^{s+t}(X/S^1). 
\]
We summarize these facts together in the Theorem below. 
\begin{prop}\label{hfpss}
Let $X$ be a free $S^1$-space, and $E$ a homotopy commutative ring spectrum. Then, there is a conditionally convergent multiplicative spectral sequence
\[
E_2^{s,t} = H^s( S^1; E^t(X)) \cong \Z[y] \otimes E^t(X) \implies E^{s+t}(X/S^1). 
\] 
1) If $E$ is complex orientable, the class $y$ is a permanent cycle. \\ 
2) The differential $d_r$ changes the grading by $(s,t)\mapsto (s+r, t-r+1)$. \\ 
3) If $X$, $X/S^1$ are finite CW complexes, and $E$ is complex orientable, the spectral sequence is strongly convergent. 
\end{prop}

\begin{proof}
The degree of the differentials follow from the construction of the exact couple for the spectral sequence. We also have the map $X_+\to S^0$ which gives a map $E \to F^{S^1}_{tr}(X_+,E)$ which is $S^1$-equivariant. Thus we have a map between the homotopy fixed point spectral sequences which maps the classes $y$ to one another, and so, 1) follows from the identification in Example \ref{z=e}. 

It remains to prove 3). For this, we show that for $k$ sufficiently large, $y^k$ lies in the image of a differential. It will then follow that for $r$ sufficiently large the classes $y^m$ and it's $\pi_\ast E$ multiples are $0$ in the $E_r$-page for $m\geq k$. Therefore, the $E_r$-page will be concentrated in the columns between $1$ and $k$, and $E_\infty=E_r$ by increasing $r$ further if necessary. Hence, the spectral sequence converges strongly. 

The space $X/S^1$ being finite dimensional implies that the classifying map $X/S^1 \to BS^1$ (for the $S^1$-bundle $X\to X/S^1$) factors through a finite skeleton. Hence, we have an equivariant map $X\to S(\C^{k+1})$ for some $k$, and thus a map $F^{S^1}_{tr}(S(\C^{k+1})_+,E) \to F^{S^1}_{tr}(X_+,E)$. As $E$ is complex orientable, 
\[
\pi_\ast F^{S^1}_{tr}(S(\C^{k+1})+,E) \cong \pi_\ast F(\C P^k, E) \cong E^{-\ast}(\C P^k) \cong \pi_\ast E [y]/(y^{k+1}) 
\]
for some choice of complex orientation $y$. Observe that the homotopy fixed point spectral sequence for the space $ES^1$ as in Example \ref{z=e} matches with the Atiyah-Hirzebruch spectral sequence for $\C P^\infty$. It follows that the class $y$ represents the complex orientation in the $E_2$-page. For $S(\C^{k+1})$ and hence also for $X$ via the equivariant map $X \to S(\C^{k+1})$, the class $y$ represents a nilpotent class whose $k+1$-power is $0$. Therefore, $y^{k+1}$ must lie in the image of a differential, and 3) follows. 
\end{proof}

\begin{exam}\label{EMhfpss}
Suppose that $E= HR$ for a commutative ring $R$, the Eilenberg-MacLane spectrum with $\pi_0 HR =R$. In this case the spectral sequence in Proposition \ref{hfpss} matches the Serre spectral sequence associated to the fibration 
\[
X \to X/S^1 \to \C P^\infty 
\]
obtained by identifying the homotopy orbits space $X_{hS^1} \simeq X/S^1$, and the classifying space $BS^1 \simeq \C P^\infty$. In this case, the spectral sequence is strongly convergent from the corresponding result for the Serre spectral sequence. Moreover, due to the fact that $X/S^1$ is a finite complex, the $E_\infty$-page vanishes beyond the dimension of $X/S^1$. 
\end{exam}
\vspace*{0.2 cm}

Next we provide a method to compute the differentials in the spectral sequence of Proposition \ref{hfpss}. In the tower of fibrations used to construct the spectral sequence, the spectrum at the bottom of the tower is $F^{S^1}_{tr}(S^1\times X_+,E)^{S^1} \simeq F(X_+,E)$. Let $Q$ denote the homotopy cofibre of the map $X \to X/S^1$. In the category of spectra, $F(\Sigma^{-1} Q_+, E)\simeq \Sigma F(Q_+,E)$ is the homotopy cofibre of the map $F(X/S^1_+, E)\to F(X_+,E)$. In view of the commutative square  
\[
\xymatrix{ F(X/S^1_+, E) \ar[r]^-{\simeq}  & F^{S^1}_{tr}(ES^1\times X_+,E)^{S^1} \ar[r] \ar[d]  & F(X_+,E) \ar@{=}[d] \\ 
& F^{S^1}_{tr}(S(\C^{k+1})\times X_+, E)^{S^1} \ar[r]   &  F(X_+,E),   
}
\]
we obtain coherent maps $\Sigma F(Q_+, E) \to \Sigma F^{S^1}_{tr}(Q(k)_+, E)^{S^1}$, where $Q(k)= [S(\C^{k+1})/S^1] \wedge X_+$ is the $S^1$-equivariant homotopy cofibre of $X\times S^1 \to X\times S(\C^{k+1})$. Projecting onto the first factor gives a map $Q(k)\to S(\C^{k+1})/S^1$ which gives a map 
\[
F(\C P^k, E) \simeq F^{S^1}_{tr}(S(\C^{k+1})/S^1, E)^{S^1} \to F^{S^1}_{tr}(Q(k)_+, E)^{S^1}. 
\]
Let an element $x \in E^n(X)$ be represented by the map $S^{-n} \stackrel{x}{\to} F(X_+,E)$. Our hypothesis about such an $x$ is a factorization in the following commutative diagram for $0\leq k \leq \infty$. 
\begin{myeq} \label{hypo}
\xymatrix{ S^{-n} \ar[d]^{x} \ar[r]^-{y}  &  \Sigma F(\C P^k,E)   \ar[d] \\ 
F(X_+,E)  \ar[r]  & \Sigma F^{S^1}_{tr}(Q(k)_+,E)^{S^1} .
}
\end{myeq}
Before applying this hypothesis we note 
\begin{prop} \label{diff-zero} 
Suppose that the composite $S^{-n} \stackrel{x}{\to} F(X_+,E) \to \Sigma F^{S^1}_{tr}(Q(k)_+,E)^{S^1}$ is null-homotopic. Then, $d_r(x)=0$ for $r\leq 2k+1$. 
\end{prop}
\begin{proof}
The statement follows from the fact that the composite being null-homotopic implies that $x$ lifts in the tower of fibrations to $F^{S^1}_{tr}(S(\C^{k+1})\times X_+, E)^{S^1}$. 
\end{proof}

\begin{exam}\label{hyp-tr}
In the case $E=HR$, the spectral sequence is the Serre spectral sequence of the fibration $X \to X/S^1 \to \C P^\infty$ according to Example \ref{EMhfpss}. Note that
\[
F(\C P^k, HR) \simeq \bigvee_{1\leq i\leq k} \Sigma^{-2i} HR,
\]
so in the diagram \eqref{hypo} $y$ may be non-trivial only when $n$ is odd, and $k\geq \frac{n+1}{2}$. If $n$ is odd and \eqref{hypo} holds for $k=\frac{n+1}{2}$, the class $x$ is transgressive, and $d_{n+1}(x)$ is the image of $y$ under the composite 
\[
S^{-n} \stackrel{y}{\to} \Sigma F(\C P^{\frac{n+1}{2}}, HR) \simeq \bigvee_{1\leq i \leq \frac{n+1}{2}} \Sigma^{-2i+1} HR \to \Sigma^{-n} HR.   
\]
\end{exam}

We assume now that $E$ is connective, and that \eqref{hypo} holds for $k=\infty$. In this case we have 
\begin{prop} \label{firsdiff} 
Suppose that \eqref{hypo} holds for $k=\infty$ and that $E$ is connective. Then, $d_r(x)=0$ if $r\leq n$. Further, $d_{n+1}(x)= d^H_{n+1}(q_H(x))$, where $q_H$ is the map $E\to H\pi_0 E$, and $d^H_{n+1}$ is the $(n+1)^{th}$ differential for the spectral sequence of Proposition \ref{hfpss} for $H\pi_0 E$.  (Here we observe that the spectral sequence is one of $\pi_0 E$-modules, so this allows us to interpret the last statement.)
\end{prop} 

\begin{proof}
We observe that $E$ is connective implies that $\Sigma F(\C P^k, E)$ is $(-2k + 1)$-connective (that is, the homotopy groups are $0$ is degree $\leq -2k$). Therefore, the composite 
\[ 
S^{-n} \stackrel{y}{\to} \Sigma F(\C P^\infty,E) \to   \Sigma F(\C P^k,E) 
\]  
is trivial for degree reasons if $-n\leq -2k$. From the commutative square 
\[
\xymatrix{ S^{-n} \ar[d]^{x} \ar[r]^-{y}  & \Sigma F(\C P^\infty,E) \ar[r] \ar[d] & \Sigma F(\C P^k,E)   \ar[d] \\ 
F(X_+,E)  \ar[r]  &  \Sigma F(Q_+, E)^{S^1} \ar[r] & \Sigma F^{S^1}_{tr}(Q(k)_+,E)^{S^1} .
}
\]
we deduce that the composite in the lower row is trivial. Hence, from Proposition \ref{diff-zero} we get that $d_r(x)=0$ if $r\leq n$. 

Via the map $q_H : E \to H\pi_0 E$, we observe that \eqref{hypo} also holds when we replace $E$ by $H\pi_0 E$. Therefore, in the associated spectral sequence $d_r^H(x)=0$ if $r\leq n$ and $d_{n+1}^H(x)$ is described in the formula in Example \ref{hyp-tr}. Also we need only assume $n$ is odd, as the result is vacuously true in the other case. We fix $k=\frac{n+1}{2}$ so that $n=2k-1$. With this choice of $n$ and $k$,
 \[
 \pi_{-n}\Sigma F(\C P^k, E) \cong \pi_0 E, ~\pi_{-n} \Sigma F(\C P^{k-1}, E)=0.
 \]
 It follows that the composite of $y$ to $\Sigma F(\C P^k, E)$ lifts to $y_k: S^{-n} \to \Sigma F(\C P^k/\C P^{k-1},E)$. 

The differential $d_{n+1}$ may be described as the composite 
\[
S^{-n} \stackrel{\chi}{\to} F^{S^1}_{tr} (S(\C^k)\times X_+, E)^{S^1}  \to \Sigma F^{S^1}_{tr} ([S(\C^{k+1})/S(\C^k)]\wedge X_+, E)^{S^1} 
\]
where $\chi$ is a lift of $x$ along the map $ F^{S^1}_{tr} (S(\C^k)\times X_+, E)^{S^1}  \to F(X_+,E)$. We expand this in the diagram below
{\small
\[
\xymatrix{ S^{-n} \ar[rdd]_(0.7){x} \ar[rd]^-{\chi} \ar[rrd]^(0.6){d_{n+1}(x)} &   \\ 
 F^{S^1}_{tr}(S(\C^{k+1})\times X_+, E)^{S^1} \ar[r] \ar@{=}[d] &  F^{S^1}_{tr}(S(\C^k)\times X_+, E)^{S^1} \ar[r] \ar[d]  & \Sigma F^{S^1}_{tr} ([S(\C^{k+1})/S(\C^k)]\wedge X_+, E)^{S^1}   \ar@{-->}[d] \\ 
F^{S^1}_{tr}(S(\C^{k+1}) \times X_+, E)^{S^1} \ar[r]       & F(X_+,E)  \ar[r]    &  \Sigma F^{S^1}_{tr}(Q(k)_+, E)^{S^1} .
}
\]
}
Observe that the dotted arrow in the diagram above is defined using the fact that the rows are homotopy cofibration sequences. We further compute $d_{n+1}(x)$ via the following commutative diagram
{\small
\[
\xymatrix{ S^{-n} \ar[d]^{\chi} \ar[r]^-{y_k} \ar[rd]^{d_{n+1}(x)}  & \Sigma F(\C P^k/\C P^{k-1},E) \ar[r] \ar[d] & \Sigma F(\C P^k,E)   \ar[d] \\ 
F^{S^1}_{tr}(S(\C^k)\times X_+, E)^{S^1}  \ar[r]  &  \Sigma F^{S^1}_{tr} ([S(\C^{k+1})/S(\C^k)]\wedge X_+, E)^{S^1}  \ar[r] & \Sigma F^{S^1}_{tr}(Q(k)_+,E)^{S^1} .
}
\]
}
The middle vertical map is the one which quotients out the factor $X$. Under the identification $\Sigma F(\C P^k/\C P^{k-1},E) \simeq \Sigma^{-2k+1} E$, and $\pi_{-n} \Sigma^{-2k+1} E \cong \pi_0 E$, we identify $y_k$ with $d_{n+1}^H(x)$. 
\end{proof} 

\section{The cohomology of $W_{n,k}$} \label{cohwnk}
In this section, we calculate the generalized cohomology of $W_{n,k}$ with respect to a complex oriented spectrum $E$. Later in the section, we specialize to $E=BP$, the spectrum for Brown-Peterson cohomology. 

Recall that a complex orientation for a homotopy commutative ring spectrum $E$ is a class $x\in \tilde{E}^2(\C P^\infty)$, which restricts to a generator of the free rank one  $\pi_0 E$-module $\tilde{E}^2(S^2)\cong E^0(pt)$. For a complex oriented spectrum $E$, we have \cite{Ada95} 
\[
E^\ast (\C P^n) \cong E^\ast(pt) [x]/(x^{n+1}),~~ E^\ast(\C P^\infty) \cong E^\ast(pt)[[x]]. 
\]

For the complex Stiefel manifold, the classical computations of their cohomology \cite{MiTo78} proceeds using the Serre spectral sequence as for other homogeneous spaces. It is proved that the cohomology of $W_{n,k}$ is an exterior algebra with generators in degrees $2n-2k+1, \cdots, 2n-1$. The Stiefel manifold also has a filtration 
\[
\xymatrix{ W_{n-k+1,1} \ar@{=}[d] \ar@{^{(}->}[r] & W_{n-k+3,2}  \ar@{^{(}->}[r] & \cdots \ar@{^{(}->}[r] W_{n-1,k-1}  \ar@{^{(}->}[r] & W_{n,k},  \\
S^{2n-2k+1}}
\]
where the inclusion $W_{n-1,k-1} \hookrightarrow W_{n,k}$ is given by adding the last vector $e_n$. The filtration quotients are computed using the following homotopy pushout \cite[Chapter IV]{Whi78}
\begin{myeq} \label{push-St}
\xymatrixcolsep{5pc}\xymatrix{
\Sigma(\mathbb{C}P^{n-2}_+ )\times W_{n-1,k-1}\ar@{^{(}->}[d] \ar[r]^ -\mu &W_{n-1,k-1}\ar@{^{(}->}[d]\\
\Sigma(\mathbb{C}P^{n-1}_+) \times W_{n-1,k-1} \ar[r]^-\mu &W_{n,k}.}
\end{myeq}
In order to construct $\mu$ one defines 
\[
S^1\times \C P^{n-1}\to U(n) 
\]
by $(z,L) \mapsto A(z,L)$, where $A(z,L):\C^n \to \C^n$ is the unitary transformation which multiplies the elements of $L$ by $z$ and fixes the orthogonal complement. The map $\mu$ is induced by matrix multiplication in $U(n)$ and the left action on $W_{n,k}$. From the construction of $\mu$ and the fact that $W_{n,k}\cong U(n)/U(n-k)$, one obtains the induced map $\mu_{n,k}: \Sigma \frac{\C P^{n-1}}{\C P^{n-k-1}} \to W_{n,k}$. It follows from \eqref{push-St} that 
\[
W_{n,k}/W_{n-1,k-1} \simeq \Sigma \frac{\C P^{n-1}}{\C P^{n-2}} \wedge {W_{n-1,k-1}}_+\simeq \Sigma^{2n-1} ({W_{n-1,k-1}}_+) . 
\]
In the case of ordinary cohomology, the exterior algebra generators for the cohomology of $W_{n,k}$ pullback under $\mu_{n,k}$ to $\Sigma x^{i-1}$ ($\Sigma : H^\ast X \to H^\ast \Sigma X$ is the suspension isomorphism). We now use this filtration to prove analogous results for the $E$-cohomology of $W_{n,k}$. 
\begin{prop}\label{E-coh-St}
Let $E$ be a complex oriented cohomology theory, such that there is no $2$-torsion in $E^\ast(pt)$. Then, 
\[
E^\ast(W_{n,k}) \cong \Lambda_{E^\ast (pt)}(z_{n-k+1}, \cdots, z_n),
\]
 is an exterior algebra with $|z_i| = 2i-1$. These generators satisfy \\
 1) The inclusion $W_{n-1,k-1} \hookrightarrow W_{n,k}$ sends $z_i$ to $z_i$ if $n-k+1 \leq i \leq n-1$ and sends $z_n$ to $0$. \\ 
 2) $\mu_{n,k}^\ast (z_i)=\Sigma x^{i-1}$. 
\end{prop}

\begin{proof}
We prove the results by induction on $k$, constructing the generators $z_i$ along the way. For $k=1$, the Stiefel manifold is the sphere $S^{2n-1}$, and in this case, we know that the $E$-cohomology is the exterior algebra on one generator. This starts the induction. 

In the induction step, we know that $E^\ast(W_{n-1,k-1})$ is as described in this Proposition, and attempt to derive the same for $E^\ast(W_{n,k})$ via the pushout \eqref{push-St}. This gives us the following maps between long exact sequences corresponding to the columns of \eqref{push-St}. 
{\footnotesize
\begin{myeq} \label{le-push-St}
\xymatrix{
\ar[r]^-0&\widetilde{E}^r(\Sigma^{2n-1}({W_{n-1,k-1}}_+))\ar[d]^-{id}\ar[r]^-{j^*}&E^r(W_{n,k}) \ar[d]^-{\mu^*}\ar@{->>}[r]^-{i^*}&E^r(W_{n-1,k-1})\ar[d]\ar[r]^-0&\\
\ar[r]^-0  &  \widetilde{E}^r(\Sigma^{2n-1}({W_{n-1,k-1}}_+))  \ar[r] &  E^r(\Sigma(\mathbb{C}P^{n-1}_+)\times W_{n-1,k-1})\ar@{->>}[r]    &   E^r(\Sigma(\mathbb{C}P^{n-2}_+)\times W_{n-1,k-1})\ar[r]^-0&} 
\end{myeq}
}
We now justify the various identifications described in \eqref{le-push-St}. The fact that $E$ is complex oriented implies $E^\ast(\C P^{n-1}) \to E^\ast (\C P^{n-2})$ is surjective, and the induction hypothesis gives us that $E^\ast(W_{n-1,k-1})$ is a free $E^\ast(pt)$-module. This implies that 
\[
E^r(\Sigma(\C P^{n-1}_+)\times W_{n-1,k-1}) \to E^r(\Sigma(\C P^{n-2}_+)\times W_{n-1,k-1})
\]
is surjective. This implies the identifications in the bottow row of \eqref{le-push-St}. The identifications on the top follow from the ones of the bottom row, and the fact that $E^r (W_{n-1,k-1}) \to  E^r(\Sigma(\C P^{n-2}_+)\times W_{n-1,k-1})$ is injective. It follows that we have short exact sequences 
\[
0 \to E^\ast (\Sigma^{2n-1} {W_{n-1, k-1}}_+) \stackrel{j^\ast}{\to} E^\ast(W_{n,k}) \stackrel{i^\ast}{\to} E^\ast(W_{n-1,k-1}) \to 0.
\]  
For $n-k+1 \leq i \leq n-1$, we choose $z_i \in E^{2i-1}(W_{n,k})$ so that they map to $z_i$ under $i^\ast$. The class $z_n$ is chosen so that it maps to $\Sigma x^{n-1}$ under $\mu_{n,n-1}$. From \eqref{le-push-St}, it follows that $z_n$ is a generator for the ideal of $E^\ast(W_{n,k})$ given by image of $j^\ast$. By the construction 1) and 2) follow. As $z_i$ are in odd degree and $E^\ast(pt)$ has no $2$-torsion, we have $z_i^2=0$, and \eqref{le-push-St} implies that $E^\ast(W_{n,k})$ additively matches with the exterior algebra on the $z_i$. The result now follows by induction on $k$.
\end{proof}

We now proceed to define the generators of the exterior algebra $E^\ast(W_{n,k})$ in a strict fashion which will satisfy 1) and 2) of Proposition \ref{E-coh-St}. From the proof, we note that for any classes $z_i$ satisfying 2), $E^\ast(W_{n,k}) \cong \Lambda_{E^\ast(pt)}(z_{n-k+1},\cdots, z_n)$.  Although the results in the following will have analogous consequences for any complex oriented $E$, we fix our attention to the case $E=BP$, which will be used in the following sections. Recall \cite{Rav86} 
\[
BP^\ast(pt) \cong \Z_{(p)}[v_1,v_2,\cdots],
\]
where $v_i$ denotes the Araki generators \cite[A2.2.2]{Rav86} that lie in degree $-2(p^i-1)$ (note here we are using the cohomological grading which is negative of the usual homotopy grading). We also fix from now on $x\in \widetilde{BP}^2(\C P^\infty)$ to denote the fixed orientation for a  $p$-typical formal group law over $BP^\ast(pt)$. We also assume that $x$ is such that it maps to the first Chern class under the map $\lambda: BP\to H\Z_{(p)}$.

The method of choosing the generators $y_j$ for $BP^\ast (W_{n,k})$ is by relating them to the $BP$-Chern classes $c_j^{BP}$ \cite{CoFl66}. We start with the case $k=n$, when $W_{n,n}= U(n)$. Recall that $H^*(U(n);\Z_{(p)})=\Lambda_{\Z_{(p)}}(y^H_1,\cdots , y^H_n)$  with $|y^H_j|=2j-1$, and in Serre spectral sequence for the fibration 
\[
U(n) \to EU(n) \to BU(n),
\]
$y^H_j$ transgresses to $j^{th}$-Chern class $c_j^H$. We also know that $A^*(y_j^H)=\Sigma x_H^{j-1}$, where $A:\Sigma (\C P^{n-1}_+) \to U(n)$ is induced from $(z,L)\mapsto A(z,L)$, and $x_H$ is the first $H$-Chern class of the canonical line bundle over $\C P^\infty$. Write $\sigma : \Sigma U(n) \to BU(n)$ for the adjoint of the equivalence $U(n) \simeq \Omega BU(n)$, and form the composite diagram
\[
 \xymatrix{\Sigma^2(\C P^{n-1}_+)\ar[r]^-{\Sigma A}\ar@{-->}@/_1pc/[rr]_-\phi&\Sigma U(n)\ar[r]^-\sigma &BU(n). }
\]
For a cohomology theory $E$, denote by $\phi_E^\ast$ (respectively $\sigma_E^\ast$) the map induced by $\phi$ (respectively $\sigma$) on $E$-cohomology. We have $\phi^*_H(c_j)=\Sigma^2 x_H^{j-1}$ as $\sigma^*_H(c_j)= \Sigma y_j^H$. 
\begin{prop}\label{U(n)-BP-coh}
There are classes $\tau_j \in BP^{2j} (BU(n))$ of the form 
\[ \tau_n = c_n^{BP},~~ \mbox{  and } \forall~ 1\leq j \leq n, \tau_j = c_j^{BP} + \sum_{k>j} \nu_k c_k^{BP} \]
for $\nu_k \in BP^\ast(pt)$, such that 
\[
\phi_{BP}^\ast \tau_j= \Sigma^2 x^{j-1}.
\]
 The standard map $BU(n)\to BU(n+1)$ classifying the canonical bundle plus a trivial bundle sends $\tau_j$ to $\tau_j$ for $j\leq n$, and $\tau_{n+1}$ to $0$. Define $y_j^{BP} \in BP^{2j-1}(U(n))$ by the formula $\Sigma y_j^{BP} = \sigma_{BP}^\ast \tau_j$. Then, \\
1) The classes $y_1^{BP}, \cdots, y_n^{BP}$ are generators for the exterior algebra $BP^\ast(U(n))$. \\
2) $\lambda(y_j^{BP}) = y_j^H$. (that is, the classes $y_j^{BP}$ are lifts of the cohomology classes $y_j^H$ to $BP$.)
\end{prop}

\begin{proof}
We note that using Proposition \ref{E-coh-St}, it suffices to prove the statements about $\tau_j$. Consider the following commutative diagram
\begin{myeq} \label{BP-H-gen}
 \xymatrixcolsep{4pc}\xymatrix{BP^*(\Sigma^2(\C P^{n-1}_+))\ar[d]_-\lambda & BP^*(\Sigma U(n)) \ar[l] \ar[d]_-\lambda & BP^*(BU(n))\ar[l]_-{\sigma^*_{BP}} \ar@{-->} @/_2pc  /[ll]_{\phi^*_{BP}} \ar[d]^-\lambda\\
H^*(\Sigma^2 (\C P^{n-1}_+);\Z_{(p)}) & H^*(\Sigma U(n);\Z_{(p)}) \ar[l]  & H^*(BU(n);\Z_{(p)})\ar[l]^-{\sigma^*_H} \ar@{-->} @/^2pc  /[ll]^{\phi^*_H}. 
} 
\end{myeq}
We have $\lambda(c^{BP}_j)=c_j$, and also that $\lambda$ maps the complex orientation of $BP$ to that of $H$. It readily follows that $\phi^*_{BP}(c^{BP}_j) - \Sigma^2x^{j-1}$ lies in the kernel of $\lambda$, which is the ideal $(v_1,v_2,\cdots)$. As
\[
BP^\ast (\Sigma^2 (\C P^{n-1}_+)) \cong \Z_{(p)}[v_1,v_2,\cdots ]\{\Sigma^2 1, \Sigma^2 x, \cdots,  \Sigma^2 x^{n-1}\}, 
\]
the left vertical arrow of \eqref{BP-H-gen} is an isomorphism in degree $2n$. It follows that $\phi_{BP}^\ast(c_n^{BP}) = \Sigma^2 x^{n-1}$, and so, $\tau_n^{BP}=c_n^{BP}$ maps to the element of $BP^\ast(\Sigma^2 (\C P^{n-1}_+))$ required by the Proposition. 

We proceed to construct the $\tau_j$ such that $\phi_{BP}^\ast \tau_j=\Sigma^2 x^{j-1}$. Starting from $j=n$, suppose that $\tau_{j+1}$ has already been defined. We now have $\phi_{BP}^\ast (c_j^{BP})-\Sigma^2 x^{j-1} \in (v_1,v_2,\cdots )$. For degree reasons we have, 
\[
\phi_{BP}^\ast (c_j^{BP})-\Sigma^2 x^{j-1} = \sum_{k > j} \rho_k \Sigma^2 x^{k-1}=\sum_{k>j} \rho_k \phi_{BP}^\ast (\tau_k),
\]
for some $\rho_k \in (v_1,v_2,\cdots)$. Rearranging terms and substituting the formula for $\tau_k$, we obtain an equation 
\[
\Sigma^2 x^{j-1} = \phi_{BP}^\ast(c_j^{BP} + \sum_{k>j} \nu_k c_k^{BP}),
\]
so that $\tau_j=c_j^{BP} + \sum_{k>j} \nu_k c_k^{BP}$ satisfies the required criteria. We note that $\phi_{BP}^\ast$ has image in $BP^\ast (\Sigma^2 (\C P^{n-1}_+))$ which is a suspension. It follows that the decomposable elements over $BP^\ast(pt)$ map to $0$ under $\phi_{BP}^\ast$. Also the formula $\phi_{BP}^\ast(\tau_j)=\Sigma^2 x^{j-1}$ implies that $\phi_{BP}^\ast$ induces an isomorphism when restricted to the module of indecomposables. This shows that the elements $\nu_k$ are unique, and so the classes $\tau_k$ are coherently defined over $n$ as required in the Proposition. 
%
%
%
\end{proof}

We now provide a strict definition for the generators of $BP^\ast(W_{n,k})$ following Proposition \ref{U(n)-BP-coh}. Recall that there are maps 
\[
i: W_{n-1,k-1} \to W_{n,k},~~ q: W_{n,k} \to W_{n,k-1},
\]
where $i$ adds the vector $e_n$ at the end, and $q$ forgets the last vector. We have already seen in \eqref{le-push-St} that $i^\ast$ is surjective in $BP$-cohomology. We also note that $q^\ast$ is injective. For, $q^\ast$ applied to the generators of $BP^\ast(W_{n,k-1})$ as in Proposition \ref{E-coh-St} together with a generator of $BP^\ast(S^{2n-2k+1})=BP^\ast(W_{n-k+1,1})$ satisfies 2) of Proposition \ref{E-coh-St}. This provides a tuple of exterior algebra generators for $BP^\ast(W_{n,k})$. Therefore, the quotient map $\pi: U(n)\to W_{n,k}$ is injective in $BP$-cohomology. 

\begin{prop}\label{BP-coh-St}
With notations as above, $\pi^\ast$ maps $BP^\ast (W_{n,k})$ to the subalgebra of $BP^\ast(U(n))$ generated by the classes $y_{n-k+1}^{BP}, \cdots, y_n^{BP}$. 
\end{prop}

\begin{proof}
We have a diagram of fibrations 
$$\xymatrix{U(n)\ar[r]\ar[d] &W_{n,k}\ar[d] \\
EU(n)\ar[r] \ar[d] &BU(n-k)\ar[d]  \\
BU(n) \ar@{=}[r] & BU(n), 
}$$
which induces the commutative diagram 
\begin{myeq}\label{gen-St}
  \xymatrixcolsep{1.5pc}\xymatrix{
\cdots \ar[r] &BP^i(U(n))\ar[r]^-\delta &BP^{i+1}(EU(n),U(n))\ar[r] &BP^{i+1}(EU(n))\ar[r] &\cdots\\ 
\cdots \ar[r] &BP^i(W_{n,k}) \ar[u]^{\pi^\ast} \ar[r]^-\delta &BP^{i+1}(BU(n-k),W_{n,k})\ar[r] \ar[u]^{\pi^\ast} &BP^{i+1}(BU(n-k))\ar[r] \ar[u] &\cdots \\
&&\widetilde{BP}^{i+1}(BU(n)).\ar[u]^-\alpha}  
\end{myeq}
From the construction of the classes $y_j^{BP}$ we have, $\pi^\ast \alpha (\tau_j)=\delta(y_j^{BP})$. On the other hand, if $j>n-k$, the class $\alpha(\tau_j)$ maps to $0$ in $BP^\ast(BU(n-k))$. The map $BP^\ast BU(n)\to BP^\ast BU(n-k)$ is the map on $BP$-cohomology associated to the standard inclusion $BU(n-k)\to BU(n)$ classifying the sum of the canonical bundle with $k$-copies of a trivial bundle. This maps $c_j^{BP}$ to $0$ for $j>k$, and hence, from the formula in Proposition \ref{U(n)-BP-coh}, the classes $\tau_j$ to $0$ if $j>k$. It follows that for $j>k$, there are classes $y_j \in BP^\ast(W_{n,k})$ such that $\alpha(\tau_j)= \delta(y_j)$, and from \eqref{gen-St} that $\pi^\ast(y_j)= y_j^{BP}$. Also the property $\phi_{BP}^\ast(\tau_j)=\Sigma^2 x^{j-1}$ implies that the classes $y_j$ satisfies 2) of Proposition \ref{E-coh-St}. This result follows readily.  
\end{proof}

\section{$BP$-cohomology of $PW_{n,k}$} \label{cohpwnk}
In this section, we describe the $BP$-cohomology ring of $PW_{n,k}$ using the homotopy fixed point spectral sequence (Proposition \ref{hfpss}). This is a strongly convergent spectral sequence 
\begin{myeq}\label{hfpsswnk}
E_2^{s,t} = \Z[x] \otimes BP^t(W_{n,k}) \implies BP^{s+t}(PW_{n,k})
\end{myeq}
Recall that 
\[
BP^\ast(W_{n,k}) \cong \Z_{(p)}[v_1,v_2,\cdots ][y_{n-k+1}, \cdots, y_n]/(y_{n-k+1}^2, \cdots, y_n^2)
\]
by Proposition \ref{BP-coh-St}. We start with a proposition describing the initial differential on the classes $y_j$. 
 
 \begin{prop}\label{diff-gen}
 In the spectral sequence \eqref{hfpsswnk}, the differentials on $y_j$ are described by
 \[
  d_{r}(y_j)= \begin{cases} 0 & \mbox{if } r < 2j\\
                     \binom{n}{j}x^j & \mbox{if } r=2j. \end{cases}
\]
 \end{prop}

\begin{proof}
The proof will follow from the existence of a diagram as in \eqref{hypo}. We have the commutative diagram
$$ \xymatrix{W_{n,k}\ar[r]\ar[d] &PW_{n,k}\ar[r]\ar[d] &\mathbb{C}P^\infty\ar[d]^-f\\
W_{n,k}\ar[r] &BU(n-k)\ar[r] &BU(n),} $$
in which the rows are fibrations \cite{AGMP}. The map $f$ classifies the $n$-fold Whitney sum of universal canonical complex line bundle. It leads to the following diagram with commutative squares
$$ \xymatrix{\widetilde{BP}^{2j-1}(W_{n,k})\ar[r]\ar[d] &\widetilde{BP}^{2j}(BU(n-k),W_{n,k})\ar[d] &\widetilde{BP}^{2j}(BU(n))\ar[l]\ar[d]^-{f^*}\\
\widetilde{BP}^{2j-1}(W_{n,k})\ar[r] &\widetilde{BP}^{2j}(PW_{n,k},W_{n,k}) &\widetilde{BP}^{2j}(\mathbb{C}P^\infty).\ar[l]} $$

Suppose that the class $\tau_j$ of Proposition \ref{U(n)-BP-coh} is mapped to $\psi_j$ under $f^*$. In the first row, the image of $y_j$ and image of $\tau_j$ coincide (Proposition \ref{U(n)-BP-coh}), hence, the same must happen in the bottom row leading to the following homotopy commutative diagram 
\begin{myeq}\label{BP-hypo}
\xymatrix{S^{-j}\ar[r]^-{\psi_j}\ar[d]^-{y_j} &\Sigma F(\mathbb{C}P^\infty,BP)\ar[d]\\ F(W_{n,k},BP )\ar[r] &\Sigma F(PW_{n,k}/W_{n,k},BP ).}\end{myeq}
The description of $\tau_j$ in Proposition \ref{U(n)-BP-coh} leads to the following form for $\psi_j$ 
\begin{myeq}\label{psi-for}
 \psi_j = \binom{n}{j}x^j + \sum_{k>j}\nu_k\binom{n}{k}x^k .
\end{myeq}
Now apply Proposition \ref{firsdiff} to get $d_r(y_j)=0$  if $r<2j$, and $d_{2j}y_j$ is determined from the corresponding spectral sequence over $H\Z_{(p)}$. This may be computed as in \cite{AGMP} to be $d_{2j} y_j= \binom{n}{j}x^j$. Hence the result follows. 
\end{proof}

We now proceed to compute the $E_\infty$-page of the spectral sequence. The main idea here is that \eqref{BP-hypo} may be used to determine all the differentials on the classes $y_j$.
\begin{prop}\label{Einftypage}
The $E_\infty$-page of the spectral sequence \eqref{hfpsswnk} is given by
$$E_\infty= \Lambda_{BP^*(pt)}(\gamma_{n-k+2},\cdots,\gamma_n)\otimes_{BP^*(pt)} BP^*(pt)[[x]]/I  $$
where $\gamma_j$ are  certain elements in $BP^*(W_{n,k})$ with $\text{deg}(\gamma_j)=2j-1$, and $I$ is the ideal of $BP^*[[x]]$ generated by the set $\{\binom{n}{j}x^j\vert n-k<j\leq n \}$.
\end{prop}

\begin{proof}
The class $x$ is a permanent cycle by Proposition \ref{hfpss}. The multiplicative structure determines all the differentials once they are known on the classes $y_j$. We notice that $E_{2n+1}$ is the $E_\infty$-page because  $d_{2n}(y_n)= x^n $ (Proposition \ref{diff-gen}) and so all the higher powers of $x$ are killed in the $E_{2n}$-page.

From Proposition \ref{diff-gen}, we see that the first non-trivial differential is $d_{2(n-k+1)}$ and  the generator $y_{n-k+1}$ and all its multiples do not survive to the next page since 
$$d_{2(n-k+1)}(y_{n-k+1})=\binom{n}{n-k+1}x^{n-k+1}. $$  
For the $y_j$ of higher degree, it may happen that the first non-trivial differential on it 
\[ d_{2j}y_j = \binom{n}{j}x^{n-k+1}\]
may be zero. This precisely happens when $\binom{n}{j}$ lies in the ideal generated by $\binom{n}{i}$ for $n-k+1\leq i < j$ inside $\Z_{(p)}$. This condition may be interpreted in terms of $p$-adic valuations of these numbers. We then obtain a multiple $p^s y_j$ on which the differential is $0$, determined by the formula $s+v_p(\binom{n}{j})= \min_{n-k+1\leq i <j}v_p(\binom{n}{i})$. The class $p^s y_j$ may now support higher order differentials. Their formula is determined by computing $p^s\psi_j$ using \eqref{BP-hypo} in the form of \eqref{hypo}
\[
\xymatrix{S^{-j}\ar[r]^-{p^s\psi_j}\ar[d]^-{p^sy_j} &\Sigma F(\C P^N,BP)\ar[d]\\ 
F(W_{n,k},BP )\ar[r] &\Sigma  F^{S^1}_{tr}([S(\C^{N+1})/S^1] \wedge {W_{n,k}}_+,BP)}
\]
for $N> 2j$. According to the formula \eqref{psi-for}, the next possible differential is 
\[
d_{2j+2}(p^s y_j)= p^s\nu_{j+1} \binom{n}{j+1}x^{j+1} = p^s \nu_{j+1} d_{2j+2}(y_{j+1}).
\]
We now rectify this class as $p^sy_j - p^s \nu_{j+1} y_{j+1}$ and obtain a cycle. This process continues until we reach the $E_{2n+1}$-page following which there are no further non-zero differentials. 

We now formalize the above process by writing down a series of modifications to produce the element $\gamma_j$. Starting with $\gamma_{j}^{(2)}:= y_{j}$, in $r$-th step of the modification, the modified element will be denoted by $\gamma_{j}^{(r)}$. Below we describe  transformations, exactly one of which will be performed to produce $\gamma_{j}^{(r+1)}$ from $\gamma_{j}^{(r)}$.
\begin{enumerate}
\item  If $d_r\gamma_{j}^{(r)}=0$, then it survives to the next page and we call that element $\gamma_{j}^{(r+1)}$. 
\item  If $r=2j$ and $\gamma_{j}^{(2j)}=y_{j}$ and $d_r(\gamma_{j}^{(r)}) = \binom{n}{j} x^{j}$. Define $s$ by the formula $s+v_p(\binom{n}{j})= \min_{n-k+1\leq i <j}v_p(\binom{n}{i})$, and declare $\gamma_{j}^{(r+1)}=p^s \gamma_{j}^{(r)}$.
\item If $r>2j$, and $d_r(\gamma_j^{(r)})\neq 0$, then we know $r$ is even, and there is a $BP^\ast$-multiple of $y_{\frac{r}{2}}$ mapped by $d_r$ onto the same class (this follows from the formula for $\psi_j$ in the same way as demonstrated for $p^sy_j$ above). That is,  $d_r(\gamma_j^{(r)})= \lambda d_{r}(y_{\frac{r}{2}})$. We declare $\gamma_{j}^{(r+1)}=\gamma_{j}^{(r)}-\lambda y_{\frac{r}{2}}$.
\end{enumerate}

We finally write $\gamma_j = \gamma_j^{(2n+1)}$ which survives to the $E_\infty$-page. Hence, we have shown that the $0$-th column of the $E_\infty$-page is $\Lambda_{BP^*}(\gamma_{n-k+2},\cdots,\gamma_n)$ . Also on the $E_\infty$-page the ideal generated by  $\{\binom{n}{j}x^j\vert n-k<j\leq n \}$ goes to $0$, as each of the generators are hit by the differentials $d_{2j}(y_j)$. This completes the proof. 
\end{proof} 

It remains now to solve the additive and multiplicative extension problems to obtain $BP^\ast PW_{n,k}$ from the expression in  Proposition \ref{Einftypage}. In the following lemma, we show that the part $BP^\ast(pt)[x]/I$ forms a subalgebra of $BP^\ast PW_{n,k}$. 
Recall the fibration $ W_{n,k}\to PW_{n,k}\stackrel{p}{\to}\C P^\infty$.  We prove
\begin{lemma}\label{kercpinfty}
The kernel of the map $p^*:BP^*(\C P^\infty)\to BP^*(PW_{n,k})$ contains the ideal $I$ generated by $\{\binom{n}{j}x^j \vert n-k<j \leq n \}$ in $BP^*(pt)[[x]]$.
\end{lemma}

\begin{proof}
The proof goes by induction on $k$. For $k=1$, the fibration is up to homotopy the following sequence   
$$ \xymatrix{W_{n,1}=S^{2n-1} \ar[r] &PW_{n,1}= \C P^{n-1} \ar@{^{(}->}[r]^-{p} &\C P^\infty,}$$ 
so that the kernel of $p^*$ is ideal generated by $x^n$, satisfying the statement of the lemma. 
Suppose that the lemma is true for $PW_{n,k-1}$. To show the result for $PW_{n,k}$ , we consider the diagram
$$\xymatrix{PW_{n,k}\ar[rr]^-{T_{n,k}}\ar[rd]^-p\ar[dd]_-q &&BU(n-k)\ar@{^{(}->}[dd]\ar@{^{(}->}[rd]\\
&\C P^\infty \ar[rr]^<<<<<<<<<<<<<<f &&BU(n).\\
PW_{n,k-1}\ar[ur]^{p}\ar[rr]_-{T_{n,k-1}}&&BU(n-k+1)\ar@{^{(}->}[ur]} $$
In the above diagram, $f$ classifies the bundle $n\gamma$ where $\gamma$ is the canonical line bundle over $\C P^\infty$, and $q$ is induced by the $S^1$-equivariant projection $W_{n,k}\to W_{n,k-1}$. The three squares in the diagram are homotopy pullbacks.   Our aim is to understand the kernel of $q^*$. We see that $PW_{n,k}\to PW_{n,k-1}$ is, up to homotopy, the sphere bundle associated to the complex bundle classified by the map $T_{n,k-1}$. This is because $BU(n-k)$ is (up to homotopy) the sphere bundle of the canonical $n-k+1$-plane bundle over $BU(n-k+1)$. As $BP$ is complex oriented, we obtain a Gysin sequence 
\[
\cdots \to BP^\ast PW_{n,k-1} \stackrel{e^{BP}(T_{n,k-1})}{\longrightarrow} BP^\ast PW_{n,k-1} \stackrel{q^\ast}{\to} BP^\ast PW_{n,k} \to \cdots 
\]
 It follows that the kernel of $q^\ast$ is the ideal generated by $e^{BP}(T_{n,k-1})$ in $BP^\ast PW_{n,k-1}$. The bundle $T_{n,k-1}$ is obtained by lifting the composite $PW_{n,k-1} \stackrel{p}{\to} \C P^\infty \stackrel{n\gamma}{\to} BU(n)$ to $BU(n-k)$ so that $T_{n,k-1} + (k-1)\epsilon = np^\ast \gamma $. We readily compute $e^{BP}$ as the top $BP$-Chern class 
 \[
 e^{BP}(T_{n,k-1})=p^\ast c_{n-k+1}(n\gamma)=\binom{n}{n-k+1} x^{n-k+1}. 
\]
Therefore, $\binom{n}{n-k+1}x^{n-k+1}$ lies in the kernel of $p^\ast: BP^\ast \C P^\infty \to BP^\ast PW_{n,k}$. By the inductive formula for the kernel of $p^\ast : BP^\ast \C P^\infty \to BP^\ast PW_{n,k-1}$, the proof is now complete. 
\end{proof}

We now apply Lemma \ref{kercpinfty} and Proposition \ref{Einftypage} to complete the calculation of $BP^\ast PW_{n,k}$. 
\begin{thm}\label{projstbp}
 For every prime $p$, the $BP$-cohomology algebra of $PW_{n,k}$ is described additively by $BP^\ast(pt)$-module
 $$ BP^*(PW_{n,k}) \cong \Lambda_{BP^*(pt)}(\gamma_{n-k+2},\cdots,\gamma_n)\otimes_{BP^*(pt)} BP^*(pt)[[x]]/I $$
 where $\gamma_j$'s are of degree $2j-1$, $x$ is of degree $2$, and $I$ is the ideal generated by  $\{\binom{n}{j}x^j\vert n-k<j\leq n \}$. This isomorphism is also multiplicative if $p\neq 2$.
\end{thm}

\begin{proof}
Lemma \ref{kercpinfty} implies that $p^\ast$ induces a ring map of $BP^\ast$-modules $BP^\ast(pt)[[x]]/I \to BP^\ast(PW_{n,k})$. Choosing representatives for generators $\gamma_j$  of Proposition \ref{Einftypage} in the $E_\infty$-page we obtain a $BP^\ast(pt)$-module map $\Lambda_{BP^*(pt)}(\gamma_{n-k+2},\cdots,\gamma_n)$ to $BP^\ast PW_{n,k}$. The multiplication as a bilinear map on these factors gives a map 
\[
\Lambda_{BP^*(pt)}(\gamma_{n-k+2},\cdots,\gamma_n)\otimes_{BP^*(pt)} BP^*(pt)[[x]]/I  \to BP^\ast PW_{n,k}
\]
of $BP^\ast(pt)$-modules. This is an isomorphism by Proposition \ref{Einftypage} and the multiplicative structure of the spectral sequence \eqref{hfpsswnk}. Further if $p\neq 2$, we have $\gamma_j^2=0$ as $\gamma_j$ lies in odd degree. Therefore, the isomorphism is also multiplicative.
\end{proof}

We observe that in Theorem \ref{projstbp}, we do not expect the isomorphism to be multiplicative when $p=2$, as it does not even hold over $H\Z/p$ (\cite{AGMP}).
 


\section{Equivariant maps between Stiefel manifolds} \label{eqstmap}
In this section, we demonstrate how the computations of $BP$-cohomology operations may be used to rule out $S^1$-equivariant maps between the Stiefel manifolds. The results of \cite{PP13} can be improved in this way. 

\begin{mysubsection}{Applications using Steenrod operations} We start with an example using Steenrod operations in $\Z/2$-cohomology. 
The Steenrod operations on $H^\ast(PW_{n,k};\Z/2)$ are described in \cite[Theorem 1.2]{AGMP}. We have from \cite{PP13} that if there is an $S^1$-equivarint map from $W_{n,k}$ to $W_{m,l}$ with $n-k=m-l$, then 
$$ \binom{n}{n-k+1}\text{ divides } \binom{m}{m-l+1},$$
which is then used to rule out such equivariant maps in many cases when $n-k=m-l$ and $n>m$ \cite[Theorem 3.10]{PP13}. 
 The Steenrod operations allow us to rule out equivariant maps for cases where the above divisibility is valid. An example is given in the Theorem below.
\begin{thm} \label{St-op-equiv}
Suppose $r \equiv$ $-1$, $-2$, or $3$ $\pmod{9}$ and $r \equiv$ $2$, $1$, or $-2$ $\pmod{7}$, and $m=16 r -2$. Then, there is no $S^1$-equivariant map from $W_{m-3,7}$ to $W_{m,10}$.
\end{thm}

\begin{proof}
Write $n=m-3$, $k=7$ and $l=10$. Observe that the following are satisfied by these integers 
\begin{enumerate}
\item $m,l$ even and $n,k$ odd, and $m-l=n-k$. 
\item $2$ divides both $\binom{m}{m-l+1} $, $\binom{n}{n-k+1}$ but $4$ does not divide either. 
\item $\binom{n}{n-k+1} \vert \binom{m}{m-l+1}$.
\end{enumerate}
An $S^1$-equivariant map $f$ from $W_{n,k}$ to $W_{m,l}$ induces a map of fibration sequences 
$$\xymatrix{W_{n,k}\ar[r]\ar[d]_f &PW_{n,k}\ar[r]\ar[d] &\C P^\infty \ar[d] \\
W_{m,l}\ar[r] &PW_{m,l}\ar[r] &\C P^\infty. }  $$
We compare the associated Serre spectral sequences with $\Z$-coefficients in the case $n-k=m-l$.   The condition (2) implies that $f^*(y_{m-l+1})=cy_{n-k+1}$ , where $c$ is odd. This is because in those spectral sequences $y_j$ transgresses to $\binom{n}{j}x^j$ and $\binom{m}{j}x^j$ respectively. The classes $y_{n-k+1}$ and $y_{m-l+1}$ also survive in the $\Z/2$-cohomology spectral sequence by (2), and we have $f^*(y_{m-l+1})=y_{n-k+1}$. \cite[Theorem 1.2]{AGMP} implies 
\[
Sq^2(y_{m-l+1})= (m-l)y_{m-l+2}+mxy_{m-l+1} = 0,
\]
and 
\[
Sq^2(y_{n-k+1})= (n-k)y_{n-k+2}+nxy_{n-k+1}= xy_{n-k+1}.
\]
This is a contradiction. 
%
%
%
%
%
\end{proof}
\end{mysubsection}

\begin{mysubsection}{Results using $BP$-operations}
We have seen how Steenrod squares yield some results on non-existence of $S^1$-equivariant maps between complex Stielfel manifolds. We now derive stronger results  using $BP$-theory and cohomology operations associated to it. The operations we use here are the Adams operations defined via \cite[2.4]{Ara75}. These are multiplicative, stable operations with the formula 
\begin{myeq}\label{adopform}
\Psi^a_{BP}(x)=a^{-1}[a]_{BP} (x),
\end{myeq}
where $a\in \Z_{(p)}^\times$, and $[a]_{BP}$ denotes the $a$-series using the $BP$-formal group law. These operations act on the coefficient ring via $\Psi^a_{BP}(v_i)=a^{p^i-1}v_i$. 

Denote the ideal $(v_1,v_2,\cdots)$ in $BP^*(pt)=\Z_{(p)}[v_1,v_2,\cdots]$ by $J$. We fix the $\{v_i \mid i \geq 1 \} $ to be the Araki generators \cite[A2.2.2]{Rav86}. The formal group law $\mu_{BP}$ , associated to $BP$ with respect to our chosen orientation is strictly isomorphic to the additive formal group law over $BP^*(pt)\otimes \Q$ and the isomorphism is given by $BP\text{-log}$ series. The choice of generators imply that  the $BP\text{-log}$ series has the form 
$$ \text{log}_{BP}(x)=x+ \sum_{i\geq1}l_i x^{p^i}, $$ 
where $l_i$ are determined by the relations 
$$ pl_n=\sum_{0\leq i\leq n}l_iv_{n-i}^{p^i} $$
with $l_0=1$ and $v_0=p$. This implies the formula 
$$l_n=\frac{v_n}{p-p^n} \pmod{J^2}.$$ 

Now we consider the expression of \cite[Part II, Proposition 7.5]{Ada95} $\pmod{J^2}$ to obtain the following relation for the $\text{exp}_{BP}$-series  
$$ \text{exp}_{BP}(x)= x- \sum_{i\geq 1}l_ix^{p^i} = x- \sum_{i\geq 1} \frac{v_i}{p-p^i}x^{p^i} \pmod{J^2}.$$ 
This implies
\begin{myeq}\label{BP-FGL} 
\begin{split}
x+_{BP}y  &= \text{exp}_{BP}(\text{log}_{BP}x + \text{log}_{BP}y)\\
                &= \text{log}_{BP}x + \text{log}_{BP}y - \sum_{i \geq 1}l_i(\text{log}_{BP} x + \text{log}_{BP} y)^{p^i} \pmod{J^2 }\\
               &= x+y +\sum_{i\geq 1}l_i(x^{p^i}+y^{p^i} ) - \sum_{i\geq 1}l_i( x+y +\sum_{i\geq 1}l_i(x^{p^i}+y^{p^i} ))^{p^i} \pmod{J^2 }\\
               &=x+y+\sum_{i\geq 1}l_i(x^{p^i}+y^{p^i}-(x+y)^{p^i}) \pmod{J^2 }\\
               &= x+y+\sum_{i\geq 1}\frac{v_i}{p-p^i}(x^{p^i}+y^{p^i}-(x+y)^{p^i}) \pmod{J^2},
\end{split}
\end{myeq}
where by $+_{BP}$ we mean the formal sum under the formal group law $\mu_{BP}$. 
%

We now restrict our attention to $p=2$, and obtain the following reduction for $\Psi^3_{BP}$ \eqref{adopform}.
\begin{align*}
\Psi^3_{BP}(x)&=\frac{1}{3}[3]_{BP}(x)\\
&=\frac{1}{3}(x+_{BP}[2]_{BP}(x) )\\
&=\frac{1}{3}(x+_{BP} (2x+_{BP} v_1x^2+_{BP}v_2x^4+_{BP}\cdots+_{BP}v_i x^{2^i}+_{BP}\cdots ))\\
&=\frac{1}{3}(x+_{BP}2x+_{BP}(v_1x^2+\cdots+v_ix^{2^i}+\cdots)) \pmod{J^2} \\
&= \frac{1}{3}((3x+\sum_{i\geq 1} \frac{v_i}{2-2^i}(x^{2^i} + (2x)^{2^i}-(3x)^{2^i}))+_{BP}(v_1x^2+\cdots+v_ix^{2^i}+\cdots)) \pmod{J^2} \\
&= \frac{1}{3}(3x+\sum_{i\geq 1} \frac{v_i}{2-2^i}(x^{2^i} + (2x)^{2^i}-(3x)^{2^i})+(v_1x^2+\cdots+v_ix^{2^i}+\cdots)) \pmod{J^2} \\
&=x+\sum_{i\geq 1}\frac{1-3^{2^i-1}}{2(1-2^{2^i-1})}v_ix^{2^i} \pmod{J^2} \text{     [applying (\ref{BP-FGL}) ]}
\end{align*}
We note that $\frac{1-3^{2^i-1}}{2(1-2^{2^i-1})}= \alpha_i$ lies in $\Z_{(2)}^\times$, and in this notation we have 
 \begin{myeq}\label{Ad-BP}
\Psi^3_{BP}(x)= x+\sum_{i \geq 1} \alpha_i v_i
x^{2^i} \pmod{J^2}
\end{myeq}

We shall now determine the action of $\Psi^3_{BP}$ on $BP^*(W_{n,k})=\Lambda_{BP^*}(y_{n-k+1},\cdots,y_n)$ modulo the ideal $I^2$, where $I$ is the ideal of $BP^*(W_{n,k})$ generated by $y_{n-k+1},\cdots,y_n$ . 
\begin{prop}\label{Ad-gen}
$$\Psi^3_{BP}(y_j)= y_j+(j-1)\sum_{i \geq 1,~ 2^i+j-1\leq n}\alpha_i v_i y_{2^i+j-1} \pmod{I^2+J^2}$$
\end{prop}
\begin{proof}
Recall the map $\mu_{n,k}:\Sigma\frac{\C P^{n-1}}{\C P^{n-k+1}}= \Sigma P_{n,k} \longrightarrow W_{n,k}$, for which we had $\mu_{n,k}^*(y_j)=\Sigma x^{j-1}$. Hence this will give us the isomorphism 
$$ \mu_{n,k}^*: BP^*(W_{n,k})/ I^2 \longrightarrow \Sigma BP^*(P_{n,k})$$ 
By naturality of the Adams operations, $\Psi^3_{BP}$ commutes with $\mu^*_{n,k}$. The action of $\Psi^3_{BP}$ on $y_j$ is determined up to $I^2$ from the computation for $\Sigma x^{j-1}$. The Adams operation being stable, commutes with the suspension, so it is enough to compute the action of $\Psi^3_{BP}$ on $x^{j-1}$, which comes from the multiplicative structure and the formulas above.
\begin{align*}
\Psi^3_{BP}(x^{j-1}) &= (\Psi^3_{BP}(x))^{j-1}\\
&= (x+\sum_{i \geq 1} \alpha_i v_i x^{2^i})^{j-1} \pmod{J^2}  \text{    (using \ref{Ad-BP})} \\
&= x^{j-1}+ \sum_{i \geq 1} (j-1) \alpha_i
v_i x^{2^i+j-2} \pmod{J^2}. 
\end{align*}
Hence the proposition follows.
\end{proof}

We now use the action of $BP$-Adams operations to prove new results about equivariant maps between complex Stiefel manifolds. We note from \cite{PP13} that the existence of a  $S^1$-equivariant map $W_{n,k}\to W_{m,l}$ implies that $n-k\leq m-l$. It states a number of hypothesis on $n,k,m,l$ in the case $n-k=m-l$ for which equivariant maps do not exist. Proposition \ref{St-op-equiv} proves some further results for this case. We use $BP$-operations to rule out equivariant maps in some cases where $n-k<m-l$.  
\begin{thm} \label{eqstbp}
Suppose that $m,n,l, k$ are positive integers satisfying\\ 
1) $n-k<m-l$ and there is an $s$ such that $m<2^s+m-l\leq n$. \\ 
2) $2$ divides all the binomial coefficients $\binom{n}{n-k+1},\cdots, \binom{n}{m-l}$. \\
3) $2$ does not divide $\binom{m}{m-l+1}$ and $2\nmid m-l$. \\
Then, there is no $S^1$-equivariant map from $W_{n,k}$ to $W_{m,l}$.
\end{thm}
\begin{proof}
We assume the contrary that $g: W_{n,k} \to W_{m,l}$ is an $S^1$-equivariant map. This induces a map of homotopy fixed point spectral sequences, and also a compatible map between the associated projective Stiefel manifolds.  
The formula for the differentials in the homotopy fixed point spectral sequence (Proposition \ref{diff-gen}) implies that $\binom{n}{m-l+1}$ must be odd due to the hypotheses 2) and 3), and the pullback satisfies
 \[
 g^*(y_{m-l+1})= \beta y_{m-l+1} + \sum _{j > m-l+1} p_j y_j  \pmod{I^2+J^2}.  
 \]
for some $\beta \in \Z_{(2)}^\times$ and $p_j\in BP^\ast(pt)$. For degree reasons, the second term in the above expression will be of the form 
$$ \sum_{j \geq 1,~ n \geq 2^j+m-l}k_j v_j y_{2^j+m-l}$$
where $k_j \in \Z_{(2)}$. 

Now we shall compute $\Psi^3_{BP}(g^*(y_{m-l+1})) $ and $g^*(\Psi^3_{BP}(y_{m-l+1}))$  modulo the ideal $I^2+J^2$. 
\begin{myeq}\label{1}
\begin{split}
    \Psi^3_{BP}(g^*(y_{m-l+1})) &= \Psi^3_{BP}(\beta y_{m-l+1} + \sum_{\substack {j \geq 1,~ n \geq 2^j+m-l}}k_j v_j y_{2^j+m-l}  ) \pmod{I^2+J^2} \\
    &= \beta (y_{m-l+1} +  (m-l)\sum_{i \geq 1,~2^i+m-l\leq n}\alpha_i v_i y_{2^i+m-l}) +\\& \sum_{\substack {j \geq 1,~ n \geq 2^j+m-l}}k_j \Psi^3_{BP}( v_j) \Psi^3_{BP}( y_{2^j+m-l}) \pmod {I^2+J^2} \\
    &= \beta (y_{m-l+1} +  (m-l)\sum_{i \geq 1,~n\geq 2^i+m-l}\alpha_i v_i y_{2^i+m-l}) + \\ & \sum_{\substack {j \geq 1 ,~ n \geq 2^j+m-l}}k_j \cdot 3^{2^j-1} v_j \cdot y_{2^j+m-l} \pmod {I^2+J^2}. 
\end{split}
\end{myeq}
On the other hand, we have 
\begin{myeq}\label{2}
\begin{split}
    g^*(\Psi^3_{BP}(y_{m-l+1}))&= g^*(y_{m-l+1} +  (m-l)\sum_{\substack{i \geq 1,~m \geq 2^i+m-l}}\alpha_i v_i y_{2^i+m-l}) \pmod{I^2+J^2}\\
    &= \beta y_{m-l+1} + \sum_{\substack {j \geq 1,~ n \geq 2^j+m-l}}k_j v_j y_{2^j+m-l}\\ & + \sum_{\substack{i \geq 1,~ m \geq 2^i+m-l}}\alpha_i v_i g^*( y_{2^i+m-l}) \pmod{I^2+J^2}.
\end{split}
\end{myeq}
Note that for degree reasons, 
\[
\alpha_i v_i g^*( y_{2^i+m-l})=\nu\alpha_i v_i y_{2^i+m-l} \pmod{I^2+J^2},
\]
for some $\nu\in \Z_{(2)}$. Since $\Psi^3_{BP}(g^*(y_{m-l+1}))=g^*(\Psi^3_{BP}(y_{m-l+1}))$, the coefficients for $y_{2^s+m-l}$ (for $s$ as in 1)) in the expressions \eqref{1}  and \eqref{2} must be the same modulo the ideal $I^2+J^2$. 
%
This implies
\begin{align*}
\beta (m-l)\alpha_s+3^{2^s-1}k_s&=k_s \\
\implies \beta (m-l) &= 2(1-2^{2^s-1})k_s .
\end{align*}
This contradicts the fact that $\beta (m-l) \in \Z_{(2)}^\times$. 
Hence no such $S^1$-equivariant map  $g$ can exist.
\end{proof}

\begin{exam}
One may easily figure out values of $m,n,l,k$ for which the hypothesis of Theorem \ref{eqstbp} are satisfied. For example putting $k=n$ and $m-l+1=2$, we obtain : If $n$ is even and  $\binom{m}{2}$ odd, and there is some $s$ such that $m<2^s+1 \leq n$, then, there is no $S^1$-equivariant map from $W_{n,n}$ to $W_{m,m-1}$.
\end{exam}

\end{mysubsection}

\vspace*{1cm}

\end{document}